\documentclass[Preprint,12pt]{elsarticle}

\usepackage{graphicx}
\usepackage{amsmath}
\usepackage{amssymb}
\usepackage{latexsym}
\usepackage{subfigure}
\usepackage{epstopdf}

\newcommand{\be}{\begin{equation}} \newcommand{\ee}{\end{equation}}
\newcommand{\barr}{\begin{array}} \newcommand{\earr}{\end{array}}
\newcommand{\bitem}{\begin{itemize}}
\newcommand{\eitem}{\end{itemize}}

\newcommand{\dstyle}{\displaystyle}

\newtheorem{thm}{Theorem}[section]
\newtheorem{lemma}[thm]{Lemma}
\newtheorem{alg}{Algorithm}

\begin{document}

\begin{frontmatter}

\title{A Stationary Accumulated Projection Method for Linear System of Equations}

\author[au1]{Wujian Peng\corref{cor1}}    \ead{wpeng@zqu.edu.cn}

\author[au2]{Shuhua Zhang}     \ead{shuhua55@126.com}

\address[au1]{Department of Math. and Info. Sciences,
          Zhaoqing Univ.,  Zhaoqing,  China,526061 }

\address[au2]{ Department of Math.,
           Tianjin Univ. of Finance and Economics, Tianjin, China,300204}

\cortext[cor1]{Corresponding Author}

\begin{abstract}
   It is shown in this paper that, almost all current prevalent iterative \mbox{methods} for solving linear system of equations can be classified as what we called extended Krylov subspace methods.
   In this paper a new type of iterative methods are introduced which do not depend on any Krylov subspaces. This type of methods are based on the so-called accumulated projection technique
   proposed by authors. It overcomes some shortcomings of classical Row-Projection technique and takes full advantages of the linear system. Comparing with traditional Krylov subspace methods which always depend on the matrix-vector multiplication with some fixed matrix, the newly introduced method (SAP) uses  different projection  matrices  which differ in each step in the iteration process to form an approximate solution. More importantly some particular accelerative schemes (named as MSAP1 and MSAP2) are introduced to improve the convergence of the SAP method. Numerical experiments show some surprisingly improved convergence behavior; some superior experimental behavior of MSAP methods over GMRES and block-Jacobi are demonstrated in some situations.
\end{abstract}


\end{frontmatter}

\section{Introduction}

  Linear systems of the form
  \be\label{eq:1}
  Ax =b
  \ee where $A\in R^{n \times n}$ is nonsingular arise from tremendous mathematical applications and are the
  fundamental objects of almost every computational process. From the very ancient
  Gaussian elimination to state-of-the-art methods like GMRES, PCG, Bicgstab (\cite{Axelson1,Axelson2,Vorst})
  as well as Multigrid method (\cite{Hackbusch_MG,Hackbusch_It}), numerous solvers of
  linear systems have been introduced and studied in
  extreme detail. Basically all solvers fall into two categories: Direct
  methods  and   iterative methods.

  Although some state-of-the-art direct methods can be applied to solve systems with pretty large amount of unknowns  (\cite{templates,Duff1997}) in some situations, for even larger scale sparse systems (say, with unknowns up to  a few \mbox{millions}) one can resort to the LGO-based solver (\cite{pengXueBao,pengDDM2009}) recently introduced by authors,
  iterative methods are the  only option available for many practical \mbox{problems}. For example, detailed three-dimensional multiphysics \mbox{simulations} lead to linear systems comprising hundreds of \mbox{millions} or even billions of unknowns, systems with several millions of unknowns are now routinely encountered in many applications, making the use of iterative methods virtually mandatory.

  Traditional iterative methods are classified as stationary and non-stationary methods. Stationary methods usually take the form:
      \be     \label{eq:stationaryIt}
           x^{k+1} = G x^k + v,  k = 0,1,2, \cdots.
      \ee
      where $v$ is a fixed vector and $x^0$ as the first guess.

  Excellent books covering the detailed analysis of error and convergence of these methods include works by Varga \cite{Varga} and David Young \cite{David}, etc.

  More recent iterative methods can be classified as Krylov subspace \mbox{methods} (or non-stationary methods in some literature)\cite{Saad}. Krylov subspace \mbox{methods} take the following form
   \be
       x_{k} = x_0 + y_k, k =1,2,\cdots
   \ee
where  $x^0$ is an initial guess and the correction vector $y_k$ belongs to a so-called Krylov subspace $$ K_m(G,r_0) \equiv span\{ r_0, Gr_0, G^2 r_0, \cdots, G^{m-1}r_0  \}. $$
By assuming different strategies for seeking $y_k$ from $ K_m(G,r_0)$ with $G$ usually taken as $A$ or $A^\prime$, one gets a variety of iterative methods such as CG, GMRES, BiCG, FOM, MNRES, SYMMLQ, etc.(\cite{Golub,Saad,Vorst}).

As a matter of fact, if we would refer extended Krylov subspace methods as those at each stage of iteration either the approximate solution or the correction vectors always come from  Krylov subspaces with a few (one or two) fixed generator matrices (By a ``generator" matrix to Krylov subspace $K_m(A,v)$ we mean
matrix $A$ here), then the traditional stationary iterative methods such as Jacobi, Gauss-Seidal and SOR can also be classified as extended Krylov subspace methods. Since for example one can easily see from (\ref{eq:stationaryIt}) that
$$x_{k+1}  = v + Gv + G^2v + G^3 v + \cdots + G^{k}v + G^{k+1 }x_0 \equiv y_{k+1} + z^{k+1}$$
where $y_{k+1} = v + Gv + G^2v + G^3 v + \cdots + G^{k}v \in K_{k+1}(G,v)$ and $z^{k+1}=G^{k+1}x_0 \in K_{k+2}(G,x_0)$ and $x_0$ is the initial guess to the system. It has been shown in \cite{Bramley} that another well-known type of methods, the row projection method(or ART method called in CT-related techniques, a generalization of RP can also be seen in \cite{pengLP}), can also be put into the form of (\ref{eq:stationaryIt}), and thus they still belong to the category of extended Krylov subspace methods.

Krylov subspace methods can be very effective when used in case the condition number of the coefficient matrix $A$ is relatively small. However in case of $A$ having large condition numbers,  they are not effective any more or even fail to converge.

Problems with current Krylov subspace methods lie on the fact that the successive corrections to previous approximation    come from Krylov subspace with a fixed ``generator" matrix and usually some fixed starting vector (say vector $v$ in above). If we take a look at the structure of Krylov subspace $K_m(A,v)=span \{v,Av, A^2v \cdots, A^{m-1}v\}$, we can see that the base vectors of this subspace always have the form  $A^k v$, which are increasingly closer to the subspaces formed by the eigenvectors corresponding to the eigenvalues with the largest magnitudes. For the sake of simplicity, we will refer them as generalized eigenspace denoted by $L_s(A)$, i.e,
\be \label{def:Ls}
   L_s(A) = span\{v_1,v_2, \cdots, v_s\}
\ee
 where $Av_i = \lambda_i$ and  $ |\lambda_i | \ge |\lambda_{i+1}|, $ for $i = 1,2,\cdots, n$.  It is thus
 inefficient to find a good ``approximation" to the error vector $e_m\equiv x - x_0$ in such a subspace when $e_m$
 contains rich eigenvector components corresponding to the smallest eigenvalues in magnitude.
 Especially when vector $e_m$ is almost perpendicular to the Krylov subspace $K_m(A,v)$.

It is thus always desirable for us to use some types of preconditioning when we apply Krylov subspace iterative methods to solve linear system of equations, especially for large scale computing. Though  numerous
preconditioning techniques are exploited in recent decades and
some of them turn out to be extremely efficient in some special situations, there does not exist a simple preconditioning technique
which can be applied in  general cases. Another important factor is, all preconditioning techniques can be traced back to certain algebraic iterative schemes (\cite{Benzi,Vorst,Wesseling}).

Our motivation here is to develop a set of purely
algebraic algorithms that can in someway overcome the difficulties
arising in the Krylov subspace methods.
Instead of seeking corrections from certain Krylov subspaces when solving system (\ref{eq:1}),
our new approach always tries to get a sequence of projection vectors $\{v_k\}$  of the solution $x$ with each $v_k$ guaranteed to be closer to $x$ than its predecessor $v_{k-1}$ by the so-called ``accumulated projection" technique.
More importantly We will also develop some accelerating techniques to improve the convergence of our
iterative methods.

\section{An Accumulated Projection Idea}

Let's start from a simple projection idea. Assume $x$ is the solution to (\ref{eq:1}) and we can get the projection vector $p$ of $x$ onto some subspace easily. Of course $p$ can be used as an approximation to $x$ and the error vector $e = x -p$ satisfies
$$
A e = r
$$
 where $r=b - Ap.$

To get a projection vector $p$ of $x$ onto any subspace $W$, we need the information of all inner products $ v_i'x$ where $v_i\, (i=1,\cdots, m)$ are the base vectors of $W$, i.e., $W = span\{v_1, v_2, \cdots, v_m\}$.  From system (\ref{eq:1}) we see that actually any groups of row vectors  $\{A_{i_1},A_{i_2},A_{i_3},\cdots, A_{i_m} \}$ can be used as the base vectors
where $A_{i_k}$ stands for the $i_k$th row vector of matrix $A$. To get a better approximation $x_{l+1}$ to $x$ provided that $x_{l}$ is given, people usually use the residual equation
\be\label{eq:residual}
    A e_{l} = r_{l}
\ee to find a correction vector $y$ so that $x_{l+1}=x_{l}+y$ is ``closer" to $x$ in some measurement while $y$ is obtained in the same way as $ x_{l} $ is calculated, as done in Row Projection Methods \cite{Bramley}.

Our idea here goes like the following. Instead of using the residual \mbox{equations} to get a correction vector, we use the given approximation $p_l$ (on subspace $W_l$) and use it as a base vector to form another subspace $W_{l+1}$ together with a group of different row vectors selected from matrix $A$ since $p_l'x$ can also be calculated when $p_l$ is formed. The projection vector $p_{l+1}$ of $x$ onto $W_{l+1}$ usually have a larger length and is thus ``closer" to $x$ than $p_l$, which is proved in Lemma \ref{lemma:3-1} and Lemma \ref{lemma:3-2}. This process can be repeated until $p_l$ reaches its limit position $x$ (actually $\bar{x}$, the projection
of $x$ on to $ran(A')$). The major features of this approach are: First it takes full advantage of all given information about the exact solution $x$ (the projection information of $x$ on each base vector $A_i$ ( $A_ix= b_i$)) as well as every approximating vector $x_l$ ($x'x_l $ is recorded and used in later calculations); secondly the \mbox{whole} process uses only the original data of the system to reach a sequence of approximation $x_{l}$ which approaches to $x$ steadily, i.e., the correction vector (which is implicitly calculated) does not rely on the residual equations. It is by this reason we name it as a ``stationary" method. By this way we can avoid the negative impact of successive matrix-vector multiplication between a few fixed generator matrices and some starting vector (i.e. the term $A^k v$ in Krylov subspace construction); Thirdly, there does not exist any constant iterative matrix ( like matrix $G$ in (\ref{eq:stationaryIt}) ) between adjacent approximations $x_l$ and $x_{l+1}$.

Assume $x'v_i = b_i, \, (i=1,2$) with $b_1 \neq 0$ and $ ||v_i||=1
\, (i=1,2)$, we wish to find a real number $t$ such that the
function $f(t) $  defined by \be \label{def:f}
     f(t) = \frac{|x'v|}{||v||}
\ee
is maximized among all possible vectors in the form $v = v_1 + tv_2$. As a matter of fact, this optimization problem is equivalent to searching a vector from subspace $span\{v_1,v_2\}$ so that it is as close to $x$ as possible.

It is easy to see from analysis that the answer to the above optimization problem lies on the following conclusion.

\begin{lemma}\label{lemma:3-1}
   Let $x'v_i = b_i, \, (i=1,2$) with $|b_1| \ge |b_2|$ and $ ||v_i||=1 \, (i=1,2)$,
and $ \alpha = v_1'v_2$. Let $ s = \frac{b_2 - \alpha b_1}{ b_1 - \alpha b_2 }$. Then
      \be \label{def:fs}
          f(s) \equiv \frac{|x'(v_1+sv_2)|}{||v_1 + sv_2||}=\dstyle{\max_{t\in R}} \frac{|x'(v_1 + tv_2  )|}{||v_1 + tv_2  ||}.
      \ee
Furthermore
\be
        f(s)\ge \max\{|b_1|,|b_2|\}
     \ee
\end{lemma}
{\bf\it Proof.}  Let $$g(t) =\frac{x'(v_1 + tv_2)}{||v_1 + tv_2||}.$$ We have
  $$  g(t) = \frac{b_1 + tb_2}{\sqrt{ 1 + 2\alpha t + t^2}}.$$ Thus
   $$
   \barr{ll}
   g'(t) &= \dstyle\frac{b_2  (1 + 2\alpha t + t^2) - (b_1 + t b_2)(\alpha + t)}{  (1 + 2\alpha t + t^2)^{3/2}} \\[4mm]
     & = \dstyle\frac{b_2 - \alpha b_1 - (b_1 - \alpha b_2)t } {(1 + 2\alpha t + t^2)^{3/2}  } \\
     &=\dstyle\frac{(b_1 - \alpha b_2)(s - t) }{(1 + 2\alpha t + t^2)^{3/2}  }
  \earr
  $$
  Let $g'(t)=0$ we have the solution as $t = \dstyle\frac{b_2 - \alpha b_1}{ b_1 - \alpha b_2 }\equiv s,$ i.e., $s$ is an
  extreme point for function $f(t)$. \\
  case 1. $ b_1> \alpha b_2,$ we have $g'(t)>0$ if ($t< s$) and $g'(t)<0 $ if $t>s$. That means $g(t)$ reaches the maximal
  value at $s$. Since $g(t) \rightarrow -b_2$ when $t \rightarrow - \infty$ and $g(t) \rightarrow b_2$ when $t \rightarrow +\infty$,
  we have $g(s)\ge g(t)>-b_2$ for all $t<s$ and $ b_2 < g(t)\le g(s)$ for all $t>s$, thus function $f(t)=|g(t)|$ reaches its maximal
  value at $s$.\\
  case 2. $ b_1 < \alpha b_2,$ we have $g'(t)<0$ if ($t< s$) and $g'(t)>0 $ if $t>s$. That means $g(t)$ reaches the minimal value
  at $s$. Since $g(t) \rightarrow -b_2$ when $t \rightarrow - \infty$ and $g(t) \rightarrow b_2$ when $t \rightarrow +\infty$,
   we have $g(s)\le g(t)<-b_2$ for all $t<s$ and $ b_2 > g(t)\ge g(s)$ for all $t>s$, thus we have $f(t)=|g(t)|$ reaches its
    maximal value at $s$.\\
  Thus in both cases we have $f(s)> |b_2|$. Since $f(0) = | g(0)| = |b_1|$ and
  $f(s)$ is the maximal value of $f(t)$, thus we also have $f(s)> | b_1|$. See figure 1. \hfill $ \Box$

\begin{figure} \centering
\subfigure[case 1] { \label{fig:a}
\includegraphics[height=5cm,width=.45\columnwidth]{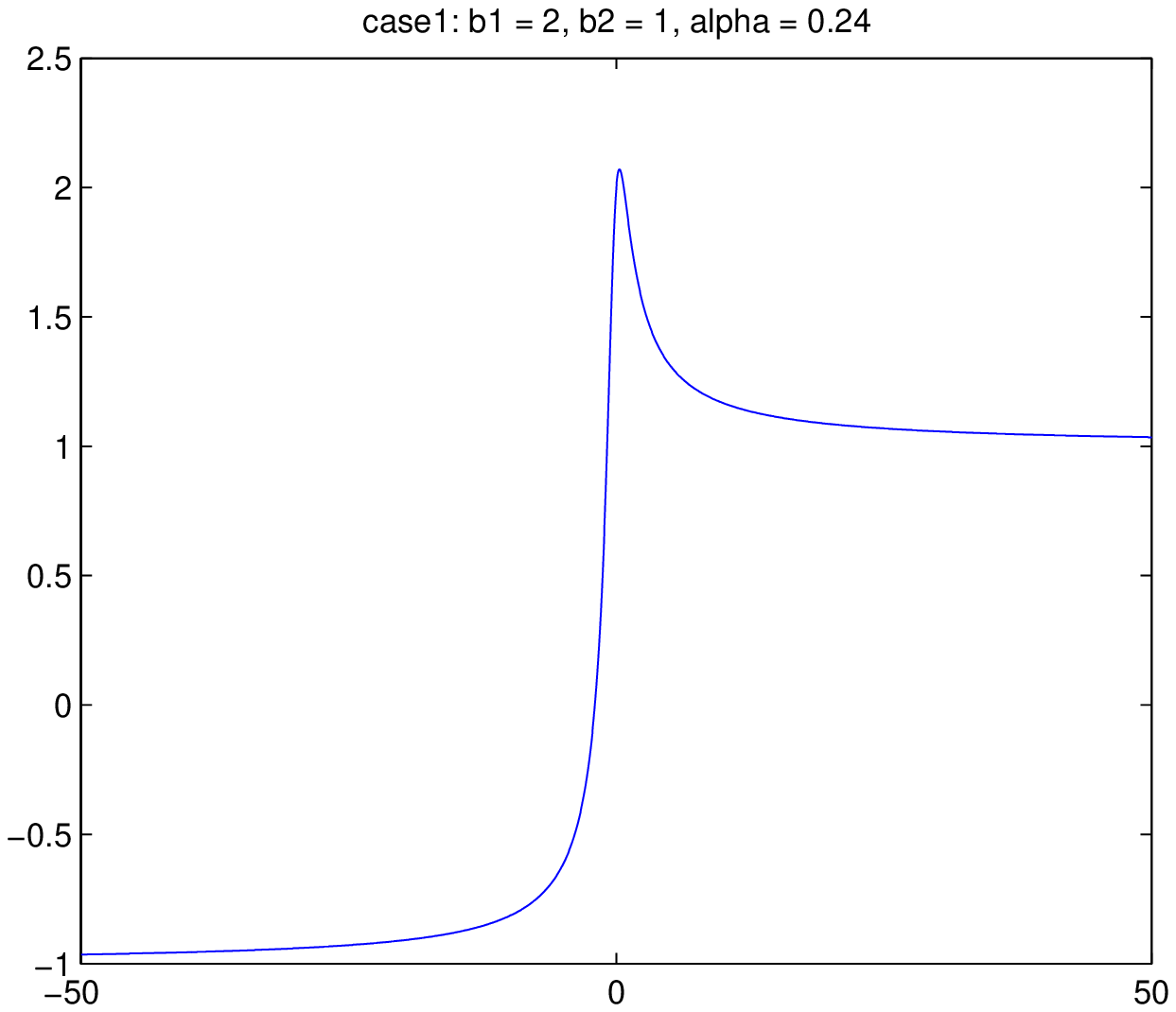}
}
\subfigure[case 2] { \label{fig:b}
\includegraphics[height=5cm,width=.45\columnwidth]{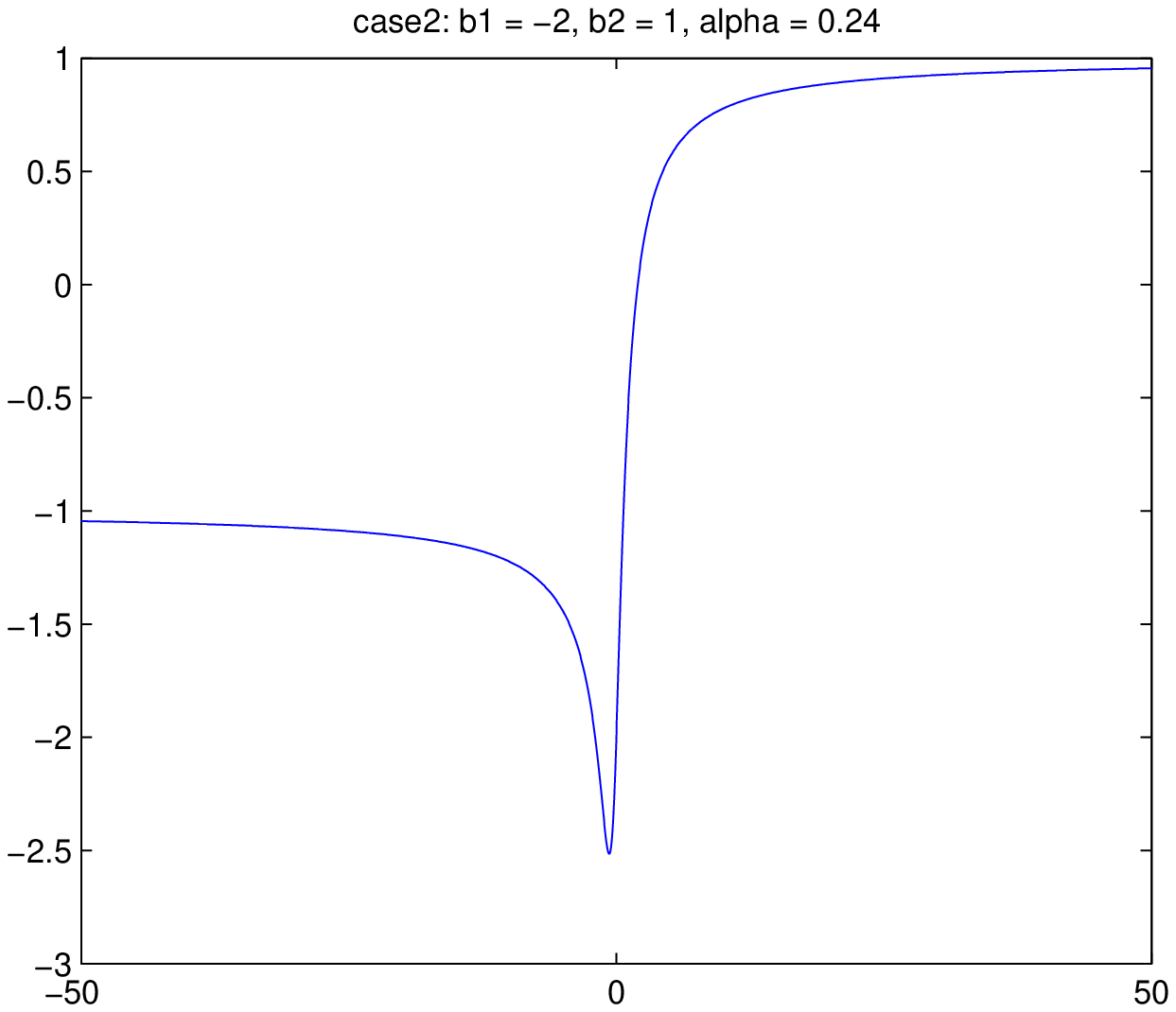}
}
\caption{ Graph of g(t) }
\label{fig}
\end{figure}

Remark: Assuming $b_1 \neq 0$, $f(s)$ can be rewritten as following (by replacing $s$ as $s = \frac{b_2 - \alpha b_1}{b_1 - \alpha b_2}$)
\be \label{fs2}
    f(s)  = \dstyle\frac{|b_1 +  s b_2|}{\sqrt{ 1 + 2\alpha s + s^2}}
          = \dstyle\frac{\sqrt{ 1- 2\alpha r + r^2}}
             {\sqrt{ 1 - \alpha^2} }|b_1|
          = |b_1|\sqrt{1 + \dstyle\frac{(r- \alpha)^2}{ 1 - \alpha^2} }
\ee where $r = b_2/b_1$.

In view of (\ref{fs2}), $f(s) \rightarrow \infty$ when $\alpha \rightarrow 1$ ( assuming that $r$ is independent of
 $ \alpha $). It is thus attempting for us to increase the length of $p_1(x)$ based on previous projection
 direction $p_0$ with $p_0'x = b_1$ by carefully selecting suitable vector $d$ with $d'x=b_2$ easily obtained
 so that $ \alpha = p_0' p_1$ is as close as possible to $1$ (i.e., the angle between $p_0$ and $p_1$ should be
  very small). However this seems to be very hard and thus we turn to an easier scheme to fulfill our task---we
  will use subspaces on which projections of $x$ are easily available. For this purpose we now generalize our
  conclusion in Lemma \ref{lemma:3-1} into the  following statement.

\begin{lemma} \label{lemma:3-2}
   Let  $ x, v_i\in R^n\, (i=1,2,\cdots, m)$,
and $ W=span\{v_1,v_2,\cdots, v_m\}$. Let $ \bar{x} $ be the projection of $x$ onto space $W$. Then
     $$ \frac{\bar{x}'x}{||\bar{x}||} = \dstyle\max_{v \in W}  \frac{|x'v|}{||v||}.$$
\end{lemma}
{\bf\it Proof}. Without loss of generality we can assume $ || x || = 1$.  By the definition of angles between vectors  we have
   $$ f(v)= \frac{|x'v|}{||v||} = \frac{|x'v|}{||v||||x||} = |\cos<x,v>|$$
   where $<x,v>$ denotes the angle between vector $x$ and $v$. Obviously $f(v)$ reaches its maximum value if and only
   if $<x,v>$ is minimized, which is true only when $v$ lies on the projection of $x$ onto subspace $W$.\hfill $\Box$

By using this result, one can always expect a searching direction $d$ on which vector $x$ has a
projection vector with larger length than any given base vectors of subspace $W = span\{v_1,v_2,\cdots, v_m\}$
with $x'v_i \, (i=1,2,\cdots, m)$ given.
Since we have $n$ vectors $A_i \,(i=1,2,3,\cdots,n)$ to form subspaces of $R^n$, this give us plenty of choices
 when it comes to construct subspaces. More importantly we can use parallel process to construct these
  subspaces and figure out projections of $x$ on each of them. Instead of using successive ``partial" projections
  which did not adequately make use of current system information, all these projections of
   $x$ can be used to construct a better approximation to the current system.

\subsection{An Accumulated Projection Algorithm}
 In this subsection we present a basic algorithm for
calculating a projection vector $p$ of $x$ to the system (\ref{eq:1}) based on current system data, i.e., the coefficient
matrix  $A$ and the right-hand side vector $b$.

In preparation, we begin with the division of all row vectors of $A$
into \mbox{groups} of vectors $\{G_i\}_1^k$, with each group $G_i$ contains
$m_i$ vectors, where $m_i \,(i=1,\cdots, k)$ are relatively small
integers satisfying $m_i<m, \, \forall  1 \le i \le k$. $m$ is a
suitable integer so that the QR factorization of matrix $A_i$ (a submatrix of $A$) formed
by all vectors in group $G_i$ is applicable; in case of sparse
coefficient matrix, QS factorization process based on LGO
method \cite{pengDDM2009} can be used and thus $m$ can be
relatively large (say, up to $O(10^5)$ or even larger). The right-hand side vector
$b$ is divided  correspondingly into vectors $b_i\, (i=1,\cdots, k)$. One
thing needs to be mentioned here is that we assume two adjacent
groups $G_i$ and $G_{i+1}$ contain about half of their vectors in
common and any row vector in $A$ must lie in at least one of the
groups, we will refer this group $\{G_i\}$ as an overlapped division
of $A$. A non-overlapped division of $A$ means the intersection of
any two groups in the division is empty.

Our approach is to  use a sequential projection process to get a
final projection vector $p$ of $x$. We begin with an initial
projection vector $p_0$ of $x$ and let it combine with all row
vectors in the first group $G_i$ to form a subspace $W_1$ of $R^n$,
and then find the projection vector $p_1$ of $x$ in $W_1$. $p_1$ is
then used to combine with all row vectors in the next group $G_2$ to
form a subspace $W_2$ so that a projection vector $p_2$ of $x$ in
$W_2$ can be obtained. The above process is repeated until all
groups are handled so that the final projection vector $p_k$ are
available.   The following algorithm gives the details.

\newpage
\begin{alg}\label{alg:AP} (An accumulated projection method-AP) The following procedure
produces an approximate vector $p$ to the solution vector $x$ which satisfies $Ax = b.$

\newcounter{apcount1}
\begin{list}{Step \arabic{apcount1}.} 
   { 
     \usecounter{apcount1} 
     \setlength{\leftmargin}{3 em}     
     \setlength{\labelwidth}{2. em}
     \setlength{\parsep}{0ex}         
     \setlength{\topsep}{1ex}         
     \setlength{\itemsep}{0.5ex}        
     \setlength{\labelsep}{.3em}     
     \setlength{\itemindent}{0em}    
     \setlength{\listparindent}{1em} 
   }
\item  Divide matrix $A$ into $k$ blocks:
 $A_1, A_2, \cdots, A_k$, divide $b$  correspondingly:  $b = b_1, b_2, \cdots,
 b_k$ ( blocks $A_i$ and $A_{i+1}$ may contain common row vectors).
 \item Initialize $p_0$ as $ p_0= \alpha A' b$ and $c_0 = \alpha ||b||^2$, where \mbox{$\alpha = ||b||^2 / ||A'b||^2$}.
\item For $i =1$ to $k$
\newcounter{apcount2}
\begin{list}{Step \arabic{apcount1}.\arabic{apcount2}} 
   { 
     \usecounter{apcount2} 
     \setlength{\leftmargin}{2 em}     
     \setlength{\labelwidth}{2. em}
     \setlength{\parsep}{0ex}         
     \setlength{\topsep}{1ex}         
     \setlength{\itemsep}{0.5ex}        
     \setlength{\labelsep}{.3em}     
     \setlength{\itemindent}{0em}    
     \setlength{\listparindent}{1em} 
   }                                                                           
   \item Construct matrix $W_i =[p_{i-1}, A_i'] $ and vector $l =[c_{i-1},b_i']'$.
   \item Compute the projection vector $p_i$ of $x$ onto subspace $ran(W_i)$  and the scalar $c_i(=x'p_i) $ as
          $$ p_i = W_i(W_i^\prime W_i)^{-1} W_i^\prime x \quad\mbox{ and } \quad c_i = l'(W_i^\prime W_i)^{-1}l.$$
   \item Go to next i.
  \end{list}
\item Output $p(=p_k)$ and $c(=c_k)$.
\end{list}
\end{alg}

It should be mentioned here that the AP algorithm depicts a successive projection process over subspace $ran(W_i) = span\{p_{i-1}, v_1,v_2, \cdots, v_{m_i}\}$ ($i=1,\cdots, k$), where $v_1, v_2, \cdots, v_{m_i}$ denotes the row vectors of submatrix $A_i$ of $A$, and $p_i$ is the projection of $x$ over subspace $ran(W_i)$ with $p_0$ stands for the initial searching direction (usually a projection vector of $x$). Obviously we have $||p_{i+1}||\ge ||p_i||$ for $ 1\le i \le k$ by Lemma \ref{lemma:3-2}.

Hence the whole AP process can be written in the matrix form as $p =  P_k x$ where $P_i$ ($i=1,\cdots ,k$) represents the projection matrix over subspace $ran(W_i$). It is easy to see that $P_i$ depends on vector $x$.  As a matter of fact, $P_k$ has the form
$$ 
   P_k = W_k(W_k'W_k)^{-1}W_k'
$$
where $W_k = [p_{k-1}, A_k']$, assuming $p_{k-1} \notin ran(A_k')$. Here and after we always use $ran(A)$ to denote the range of matrix $A$, i.e, the subspace formed by all column vectors of matrix $A$.

As a straightforward application, Algorithm \ref{alg:AP}  can be used to
solve the linear system (\ref{eq:1}) as stated in the next algorithm.\\

\begin{alg}\label{alg:SAP} (Stationary Accumulated Projection Method-SAP). Let $A \in R^{m\times n}$ , $b \in R^n$ with $m \le n$. $\epsilon$ be a given tolerance. The following procedure produces an approximation $p$ to the vector $x \in R^n$ satisfying $Ax = b$.

\newcounter{sapcount1}
\begin{list}{Step\arabic{sapcount1}.} 
   { 
     \usecounter{sapcount1} 
     \setlength{\leftmargin}{3 em}     
     \setlength{\labelwidth}{2. em}
     \setlength{\parsep}{0ex}         
     \setlength{\topsep}{1ex}         
     \setlength{\itemsep}{0.5ex}        
     \setlength{\labelsep}{.3em}     
     \setlength{\itemindent}{0em}    
     \setlength{\listparindent}{1em} 
   }

\item Initialize $s$ as $s=0$, vector $x_{0}$ as \mbox{$x_{0} = \alpha A'b$}, \mbox{$c_{0} = \alpha b'b$}, $ tol = ||b - Ax_{0}||/ ||b||$, where \mbox{$\alpha =||b||^2/||A'b||^2 . $}
\item While $ tol > \epsilon$
\newcounter{sapcount2}
\begin{list}{Step\arabic{sapcount1}.\arabic{sapcount2}} 
  {  
     \usecounter{sapcount2} 
     \setlength{\leftmargin}{3 em}     
     \setlength{\labelwidth}{2. em}
  }

        \item  Use Algorithm \ref{alg:AP} to get a projection vector $x_{s+1}$ of $x$ and $c_{s+1}(=x'x_{s+1})$ with $p$ and $c$ in step 2 replaced by  $x_{s}$ and $c_{s}$ respectively.

        \item Calculate $tol = ||b - Ax_{s}||/ ||b||;$
        \item  $s = s+1$;
       \end{list}
\item Output $p(=x_{s})$  and $c(=c_{s})$.
\end{list}
\end{alg}

{\it Remark:} In actual implementation of SAP algorithm, in order to \mbox{effectively} obtain the projection of $x$ over each subspace $ran(W_i)$ through Algorithm 1, one can store the resulting QR or LGO factors ( $Q_i$ and $R_i$ for QR, or $Q_i$ and $S_i$ for LGO respectively)  of all submatrix $A_i^\prime\, (i=1,2,\cdots, k)$ once in advance and reuse them in later projections. Although the projection matrix $P_i$ varies constantly, the projection vector can always be obtained in an economic count of flops, as it can be seen in later sections.

\begin{figure}[h]
\begin{center}
\includegraphics[height=5cm,width=10cm]{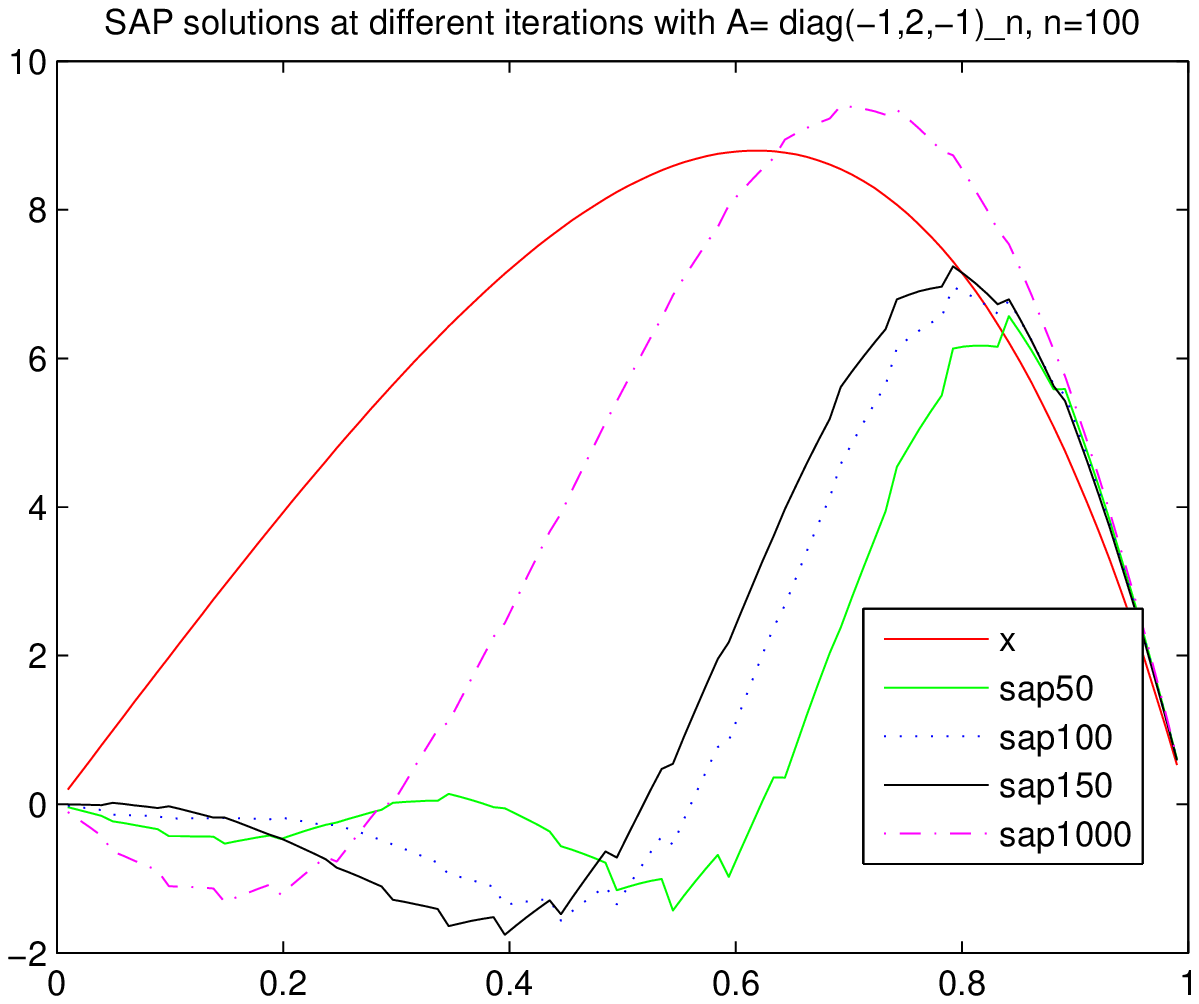 }
\caption{
Comparison of approx. solns at different \mbox{iteration} numbers
}
\label{figure:iters}
\end{center}
\end{figure}

The convergence of this algorithm is put forward to the next section. We need to point out that each sweep in step 2 is a projection process with projection matrix $P_s$ ($s = 1,2,\cdots k$) varies. Figure \ref{figure:iters} shows the comparison between approximate solutions at different iterations by this algorithm,
 and Table \ref{table:iters1} gives the
needed iteration for a convergent solution under given tolerance, where
the coefficient matrix $A$ is chosen as $A = diag(-1,2,-1)$ with $A \in R^{100\times 100}$ and the block size is chosen as $20$ when applying Algorithm \ref{alg:SAP} in this case.

\begin{table}[h]\begin{center}
\caption{SAP--iteration numbers needed for convergence}
\label{table:iters1}
\begin{tabular}{|c|c|c|c|c|c|}
\hline  tolerance  & $10^{-3}$ &$ 10^{-4}$ &$ 10^{-5} $ & $10^{-6}$ & $10^{-7}$ \\ \hline
   iter\# &    $ 724$   &     $ 872 $  &   $1020$ &   $1169 $  &   $1317$  \\ \hline
\end{tabular}
\end{center}
\end{table}

\section{Error  Analysis }

In this section we present some  analysis results for AP process and the SAP algorithms. We need to mention here that unlike classical Krylov subspace methods, the SAP method proposed here can actually be used to solve any under-determined systems.

\subsection{AP analysis}
We first present some analysis about the AP process described in Algorithm \ref{alg:AP}.
\begin{lemma} \label{lemma:AP1}

 Assume that matrix $A \in R^{m\times n}$ ($m \le n$) has full row rank, $x \in R^n$ and $b \in R^m$ where $m\le n$ satisfy $Ax = b$.  Let $A$ be divided into $k$ submatrices by its rows:  $A = (A_1',A_2',\cdots, A_k')' $ with $A_i \in R^{m_i\times n}$, and $b$ is divided as  $b = (b_1',b_2',\cdots, b_k')'$ correspondingly.   Let $\{p_i\}_1^k$ be the vector sequence produced by
 AP process (Algorithm 1).

\begin{enumerate}
      \item[(1)]    There holds for every $i=1,2,\cdots k$
           \be \label{eq:2}
            (x - p_i,p_s) = 0, \quad (s=i,i-1).
            \ee
      \item[(2)]
           Vector $ p_{i+1} - p_i$ ($i=0,1,\cdots, k-1$) is orthogonal to $p_i$, i.e.
           \be  \label{eq:3}
            (p_{i+1} -  p_i,p_i)=0
           \ee

      \item[(3)] There holds for $i=1,2,\cdots, k$
            \be  \label{eq:4}
              ||p_i||^2 + ||p_{i+1} - p_i||^2= ||p_{i+1} ||^2
            \ee

      \item[(4)]  For every $s (1 \le s \le k)$,  there holds
          \be  \label{eq:5}
            ||p_{s}||^2 = ||p_0||^2 + \sum_{i=1}^s {||p_i-p_{i-1}||^2}
          \ee
   \end{enumerate}

 \end{lemma}

{\bf\it Proof.}
  (1)  We first show that $(x-p_0,p_0) =0$. As a matter of fact, since $\alpha = b'b/(b'AA'b)$,
  we have
     $$ (x - p_0,p_0) = (x - \alpha A_1'b_1, \alpha A'b) = \alpha x'A_1'b - \alpha^2 b'A_1A_1'b
            =  \alpha b'b - \alpha b'b = 0$$
  From the fact that $p_i$ is the projection of $x$ over subspace $ran(W_i)$ with
     $W_i = [p_{i-1}, A'_i]$, for any $i ( 1\le i \le k)$ we must have
          \begin{center} $(x-p_i, p_i)=0 $ and $(x - p_i, p_{i-1})=0$
     \end{center}
      since  both $p_i$ and $p_{i-1} $ belong to $W_i$.

  (2) Note that from (\ref{eq:2}) we have
      $$ (p_{i+1} -p_i, p_i) =( ( x-p_i) - (x - p_{i+1}), p_i) =
             ( x-p_i, p_i) -(x - p_{i+1}, p_i)=0, $$
      which yields (\ref{eq:3}).

  (3) From (\ref{eq:3}) we have  $$
       \begin{array}{ll}
        || p_{i+1} - p_i||^2 &= ( p_{i+1}- p_i, p_{i+1} - p_i) \\
                              &=(p_{i+1}- p_i, p_{i+1} ) \\
                              & =  (p_{i+1}, p_{i+1}) - (p_i, p_{i+1}) \\
                              & =  (p_{i+1}, p_{i+1})  - (p_i, (p_{i+1} - p_i) + p_i)\\
                              & =    (p_{i+1}, p_{i+1}) - (p_i,p_i) \\
                              &  =  || p_{i+1} ||^2 - ||p_i||^2
       \end{array}
       $$
       from which (\ref{eq:4}) comes immediately.

  (4) Equation (\ref{eq:5}) follows from the recursive application of (\ref{eq:4}):
        $$
        \begin{array} {rl}
          ||p_s||^2 &= ||p_{s-1}||^2 + || p_s - p_{s-1}||^2 \\
                   & =||p_{s-2}||^2 +|| p_{s-1} - p_{s-2} ||^2 +   || p_s - p_{s-1}||^2 \\
                   & \cdots \\
                 & = ||p_0||^2 + ||p_1 - p_0 ||^2 + ||p_2 - p_1||^2 + \cdots +   || p_s - p_{s-1}||^2 .

        \end{array}
        $$ Proof is completed  \hfill$\Box$

Lemma \ref{lemma:AP1} actually tells the fact that the ``length" (norm) sequence $\{||p_i||\}_1^k$ of
projection vector $\{p_i\}_1^k$ actually forms a monotonically increasing sequence, and obviously $||x||$ is actually
one of its upper bounds. In order to find out how fast this sequence is increasing, we need to figure out the
detailed information of each $||p_i|| \,(i=1,2,\cdots, k)$. The following conclusion answers this question.

\begin{lemma} \label{lemma:buildpi1}

 Assume the same assumption in Lemma \ref{lemma:AP1}. Then $p_{i+1}$ has the following expression
   \be \label{eq:6}
     p_{i+1} =  \alpha_i p_i + A'_{i+1} u
   \ee
   where $u$ is
   \be\label{eq:7}
    u = \tilde{A}_{i+1} (b_{i+1} - \alpha_i A_{i+1} p_i)
    \ee
   and
   \be \label{eq:8}
    \alpha_i = \frac{x'p_i - (A_{i+1}p_i)' \tilde{A}_{i+1}  b_{i+1}    }
   { p_i'p_i - p_i'A_{i+1}'\tilde{A}_{i+1} A_{i+1} p_i }
   \ee
   and
   $$ \tilde{A}_{i+1} = (A_{i+1} A_{i+1}')^{-1}.$$
   Furthermore
    \be\label{eq:norm_pk} ||p_{i+1}||^2 = \alpha_i ||p_i||^2  + b_{i+1}' \tilde{A}_{i+1} b_{i+1}-
    \alpha_i^2(A_{i+1}p_i)'\tilde{A}_{i+1} (A_{i+1}p_i)
    \ee

 \end{lemma}
{\bf \it Proof.} \\
 It is valid to express $p_{i+1}$ in the form like (\ref{eq:6}) for some $u \in R^{m_{i+1}}$  since
 $ p_{i+1} \in ran(W_{i+1})$, where $m_i$ is the number of rows in submatrix $A_i$.

 Since $p_{i+1}$ is the projection of $x$ over subspace $W_{i+1}$, we have
    $$ A_{i+1}(x-p_{i+1})=0$$
    which leads to
    $$ b_{i+1} - \alpha_i A_{i+1} p_i - A_{i+1}A_{i+1}' u = 0$$
    from which comes (\ref{eq:7}).

 Similarly, by $ p_i'(x-p_{i+1}) =0$ we have
     $$ x'p_i - \alpha_i p_i'p_i + u'A_{i+1}p_i =0,$$
 replacing $u$ by  (\ref{eq:7}) yields (\ref{eq:8}).

 Finally from (\ref{eq:6}) we have
 \be     \label{eq:9}
     \barr{rl}
        ||p_{i+1} ||^2 & = (\alpha_i p_i + A_{i+1}' u)'(\alpha_i p_i    + A_{i+1}' u)  \\
                       & = \alpha_i^2 p_i'p_i + 2 \alpha_i p_i'A_{i+1}' u+ u'A_{i+1}A'_{i+1}u.   \\

     \earr
 \ee
 Since
 \be \label{eq:10}
 \barr{rl}
           2 \alpha_i p_i'A_{i+1}' u &= 2 \alpha_ip_i'A_{i+1}'\tilde{A}_{i+1}(b_{i+1} - \alpha_i A_{i+1}'p_i)\\
                                     & = 2 \alpha_ip_i'A_{i+1}'\tilde{A}_{i+1}b_{i+1} - 2\alpha_i^2 (A_{i+1}p_i)'\tilde{A}_{i+1}(A_{i+1}p_i)
   \earr
 \ee
 and
\be        \label{eq:11}
  \barr{rl}
            u'A_{i+1}A'_{i+1}u &= ( b_{i+1} - \alpha_i A_{i+1}p_i)'\tilde{A}_{i+1}( b_{i+1} - \alpha_i A_{i+1}p_i)\\
                               & = b_{i+1}'\tilde{A}_{i+1} b_{i+1} - 2\alpha_i b_{i+1}' \tilde{A}_{i+1} A_{i+1}p_i \\
                               & ~~+\alpha_i^2 (A_{i+1}p_i)'\tilde{A}_{i+1}(A_{i+1}p_i), \\
     \earr
 \ee
equation (\ref{eq:norm_pk}) comes from (\ref{eq:9})  (\ref{eq:10}) (\ref{eq:11}) combined. \hfill$\Box$

Lemma \ref{lemma:buildpi1} describes one way of constructing $p_{i+1}$, and detailed information about $p_{i+1}$
is revealed by (\ref{eq:norm_pk}). However a more direct approach can be used to evaluate the difference of the norms between two  consecutive projections $ p_{i+1}$ and $p_i$. These can be shown in the following
conclusion.

\begin{lemma} \label{lemma:buildpi2}
 Assume the same assumption in lemma \ref{lemma:AP1}. Let $I$ be the identity matrix in $R^n$. Then $p_{i+1}$ has the following expression
\be\label{eq:piform2}
    p_{i+1} = p_i + \bar{A}_{i+1}' v
\ee
and
\be \label{eq:pinorm}
 ||p_{i+1}||^2 - ||p_i||^2   = (b_{i+1} - (p_i'x)d)'(\bar{A}_{i+1}  \bar{A}_{i+1}')^{-1}  (b_{i+1} - (p_i'x)d)
\ee
where
   $\bar{A}_{i+1}$ is a rank-one modification of submatrix $A_{i+1}$ as
   \be\label{def:Ai_bar}
      \bar{A}_{i+1}=  A_{i+1} - d p_i' = A_{i+1}(I - u_iu_i')
   \ee
   with $u_i = p_i/||p_i||$,  $d \in R^{m_{i+1}}$ a vector taken as
   $ d =  A_{i+1}p_{i+1} /  ||p_i||^2 $
   and $v$ is defined as
   $$
      v = (\bar{A}_{i+1}\bar{A}_{i+1}')^{-1} (b_{i+1} - (x'p_i)d )
   $$
   assuming the related inverse exists.
\end{lemma}
\emph{Proof.} \\
Since $p_{i+1}$ is the projection of $x$ over subspace $W_{i+1} = ran([p_i, A_{i+1}'])$,
it can be constructed as follows.

First we modify row vectors in $A_{i+1}$ so that they are orthogonal to vector $p_i$, this
can be depicted as a rank-one modification to $A_{i+1}$ as
   $$ \bar{A}_{i+1} = A_{i+1} - d p_i',$$
where $d$ can be obtained from the fact that
  $$ \bar{A}_{i+1} p_i = 0$$
which leads to
   $$ A_{i+1} p_i - d p_i'p_i =0,$$
hence
  $$ d = A_{i+1}p_i/(p_i'p_i), $$ and
  $$\bar{A}_{i+1} = A_{i+1} - dp_i' = A_{i+1} - A_{i+1}p_ip_i'/(p_i'p_i) = A_{i+1}(I - u_iu_i'),$$
  where $u_i = p_i/||p_i||$.

Next we calculate the projection vector $\tilde{p}_{i+1}$  of $x$ over $ran(\bar{A}_{i+1})$ as
$$ \tilde{p}_{i+1} = \bar{A}_{i+1}'v,$$
where $v$ can be derived from the fact that
$$ \bar{A}_{i+1} (x - \tilde{p}_{i+1} ) = 0,
$$
which  leads to
$$
   v = (\bar{A}_{i+1}  \bar{A}_{i+1}')^{-1} (b_{i+1}- (p_i'x)d )
$$
assuming    $(\bar{A}_{i+1}  \bar{A}_{i+1}')^{-1}$ exists.

Since $\tilde{p}_{i+1} = \bar{A}_{i+1} v$ is the projection of $x$ over $ran (\bar{A}^\prime_{i+1})$ and
$\bar{A}_{i+1} p_i = 0$, we must
have $(p_i, \tilde{p}_{i+1}) =0 $. Therefore
$$
\barr{rl}
   ||p_{i+1}||^2 - ||p_i||^2 &= ||\tilde{p}_{i+1} ||^2             \\
                             & = v'\bar{A}_{i+1} \bar{A}_{i+1}'v \\
                             & = (b_{i+1} - (p_i'x)d)'(\bar{A}_{i+1}  \bar{A}_{i+1}')^{-1}  (b_{i+1} - (p_i'x)d).
\earr
$$
noting that matrix $(\bar{A}_{i+1}  \bar{A}_{i+1}')^{-1}$ is symmetric (actually positive definite symmetric).\hfill$\Box$.

\emph{Remark}: It can be shown that the length difference between $p_{i+1}$ and $p_i$ can also be written as
\be \label{diff:len}
    ||p_{i+1}||^2 - ||p_i||^2 = \tilde{x}' G \tilde{x}
\ee
where $ G =(\bar{A}_{i+1}\bar{A}_{i+1})^{-1}$ and $\tilde{x} = \bar{x} - (x'u)u$, where $\bar{x}$ denotes the
projection of $x$ on $ran(A')$ and $(x'u)u$ is the projection of $x$ (as well as $\bar{x}$) on the direction of $u= p_i/||p_i||$.

 Note that in the above lemma, we need to assume the existence of each matrix  $(\bar{A}_{i+1}  \bar{A}_{i+1}')^{-1}$. The following conclusion gives the sufficient and necessary conditions for these to hold true.

 \begin{lemma}\label{lemma:3-5}
    Let $A \in R^{m\times n} (m \le n)$ and $rank(A) = m$, $u \in R^n$ be a unit vector
    in $R^n$.  Let $\bar{A} = A(I - uu')$ and $G = \bar{A}\bar{A}'$, where $I$ denote the identity
    matrix in $R^n$. Then $G$ is nonsingular if and only if $u \notin ran(A')$.
 \end{lemma}
 \emph{Proof}.  \\
 Note that $G=\bar{A}\bar{A}'$ is invertible if and only if $ \bar{A}$ is of full row rank.

       (Necessity) Assume $G$ is invertible, we need to show that $u \notin ran(A')$. If this is not the case,
       i.e., $u \in ran(A')$, then there is a $v \in R^n$ ( $v \neq 0$) such that $ u = A'v$.
       Thus
       $$
         \bar{A}' v = (A(I- uu'))'v = (A - Auu')'u = A'v - uu'A'v = u - u(u'u) = 0
       $$
       since $||u||=1$.  This means $\bar{A}$ is not of full rank, hence $G$ is singular, a contradiction
       with our assumption.

       (Sufficiency). Assume $u\notin ran(A')$, we need to show that $G$ is invertible. As a matter of fact, if
       $G$  is not invertible, then $\bar{A}$ is not of full-row rank. Therefore there exists a nonzero vector
       $v\in R^m$
       such that $\bar{A}'v =0$. That means
        $$ 0= (A(I - uu'))'v = A'v - uu'A'v = A'v - \alpha u $$
        where $\alpha = u'(A'v)$ is a scalar. It is easy to see from here that $\alpha \neq 0$, otherwise we would have
        $A'v=0$ which means $A$ is not of full row rank. Hence $u = A'v/\alpha$, i.e., $u\in ran(A')$, this is
        contradictory    with the assumption.  \hfill$\Box$

\begin{lemma} \label{lemma:AP4}
 Assume the same assumption in Lemma \ref{lemma:AP1}. Vector sequence $p_0$,$p_1$,$\cdots$, $p_k$ are produced in
 one AP process, then
 \be \label{eq:pinorm_relation}
      ||p_i|| \le||p_{i+1}|| \quad ( i=0,1,2,\cdots, k)
 \ee
 and the equal sign holds  if and only if
   $$ A_{i+1} p_i = b_{i+1}$$
 \end{lemma}
\emph{ Proof. } Inequality (\ref{eq:pinorm_relation}) comes from (\ref{eq:4})  directly.
  We now prove the necessary condition for $||p_{i+1}|| = || p_i||$.

 \noindent (Necessity)
  Note that if $||p_{i+1}|| = ||p_i||$ holds , by (\ref{eq:4}) we must have $p_{i+1} = p_i$.
  Also from (\ref{eq:6}) we know that
    $$ p_{i+1} = \alpha_i p_i + A_{i+1}' u, $$
  thus
  \be \label{eq:tt}
     A_{i+1}' u = p_{i+1} -  \alpha_i p_i  = (1 - \alpha_i) p_i.
  \ee
  Multiplying both sides of (\ref{eq:tt}) by $A_{i+1}$ we have
  \be \label{eq:tt1}
     A_{i+1} A_{i+1}'u = (1- \alpha_i) Ap_i .
  \ee
  Note that from (\ref{eq:7}) we have
  \be \label{eq:tt2}
     A_{i+1} A_{i+1}'u =  b_{i+1} - \alpha_i Ap_i.
   \ee
  Combining (\ref{eq:tt1}) and (\ref{eq:tt2}) yields $$ Ap_i = b_{i+1}.$$

  \noindent(Sufficiency)Now we prove  $p_{i+1} = p_i$ under the assumption $ Ap_i = b_{i+1}.$

  As a matter of fact, in view of (\ref{eq:7}) and (\ref{eq:8}) we only need to show that
   $$ \alpha_i  = 1 $$
   in this case.

   Since  $ (x - p_i,p_i) = 0,$ we have
    \be\label{eq:tm1}
      x'p_i = p_i'p_i .
    \ee
    By using $ A_{i+1}p_i = b_{i+1}$ we obtain
    \be\label{eq:tm2}
       p_i'A_{i+1}' \tilde{A}_{i+1} A_{i+1}p_i = b_{i+1} \tilde{A}_{i+1} b_{i+1}
    \ee
    Hence from (\ref{eq:8}) we have
    $$
      \alpha_i = \frac{x'p_i - (A_{i+1}p_i)' \tilde{A}_{i+1} b_{i+1} }
                      { p_i'p_i -  p_i'A_{i+1}' \tilde{A}_{i+1} A_{i+1}p_i}
               = \frac{p_i'p_i - b_{i+1}' \tilde{A}_{i+1} b_{i+1} }
                      { p_i'p_i - b_{i+1}' \tilde{A}_{i+1} b_{i+1} }
                =1
                $$
    This completes the proof of the sufficient condition. \hfill $\Box$

Based on the above error analysis about AP process, a practical AP algorithm can be implemented as follows.

\begin{alg}\label{alg:APnew} (An accumulated projection method-AP) The following procedure
produces an approximate vector $p$ to the solution vector $x$ which satisfies $Ax = b.$

\newcounter{apcount1n}
\begin{list}{Step \arabic{apcount1n}.} 
   { 
     \usecounter{apcount1n} 
     \setlength{\leftmargin}{3 em}     
     \setlength{\labelwidth}{2. em}
     \setlength{\parsep}{0ex}         
     \setlength{\topsep}{1ex}         
     \setlength{\itemsep}{0.5ex}        
     \setlength{\labelsep}{.3em}     
     \setlength{\itemindent}{0em}    
     \setlength{\listparindent}{1em} 
   }
\item  Divide matrix $A$ into $k$ blocks by its row vectors:
 $A_1, A_2, \cdots, A_k$, divide $b$  correspondingly:  $b_1, b_2, \cdots,
 b_k$(Note, $A_i$ and $A_{i+1}$ may contains some common row vectors).
 \item Initialize $p_0$ as $ p_0= \alpha A' b$ and $c_0 = \alpha ||b||^2$, where \mbox{$\alpha = ||b||^2 / ||A'b||^2$}.
\item For $i =1$ to $k$
\newcounter{apcount2n}
\begin{list}{Step \arabic{apcount1n}.\arabic{apcount2n}} 
   { 
     \usecounter{apcount2n} 
     \setlength{\leftmargin}{2 em}     
     \setlength{\labelwidth}{2. em}
     \setlength{\parsep}{0ex}         
     \setlength{\topsep}{1ex}         
     \setlength{\itemsep}{0.5ex}        
     \setlength{\labelsep}{.3em}     
     \setlength{\itemindent}{0em}    
     \setlength{\listparindent}{1em} 
   }                                                                           
   \item Do QR factorization on submatrix $A_i$ : $ A_i' = Q_iR_i$.
   \item Compute vector $\tilde{b}_i = (R_i^\prime)^{-1}b_i$.
   \item Compute projection vector $x_i$ of $x$ over $ran(A_i')$:   $x_i = Q_i'\tilde{b}_i$.
   \item Store orthogonal matric $Q_i$ and vector $\tilde{b}_i$.
   \item Go to next i.
  \end{list}

\item For $i =1$ to $k$
\newcounter{apcount2nn}
\begin{list}{Step \arabic{apcount1n}.\arabic{apcount2nn}} 
   { 
     \usecounter{apcount2nn} 
     \setlength{\leftmargin}{2 em}     
     \setlength{\labelwidth}{2. em}
     \setlength{\parsep}{0ex}         
     \setlength{\topsep}{1ex}         
     \setlength{\itemsep}{0.5ex}        
     \setlength{\labelsep}{.3em}     
     \setlength{\itemindent}{0em}    
     \setlength{\listparindent}{1em} 
   }                                                                           
   \item Compute projection vector $\tilde{p}_{i-1}$ of $p_{i-1}$ on $ran(A_i')$:
   $$\tilde{p}_{i-1}= Q_i(Q_i'p_{i-1}).$$
   \item Compute vector $\bar{p}_{i-1} = p_{i-1} - \tilde{p}_{i-1}.$ 
   \item Compute the projection vector $\tilde{x}_i$ of $x$ on direction  $\bar{p}_{i-1}$:
    $$\tilde{x}_i = \beta \bar{p}_{i-1}, \mbox{ where } \beta =  (c_{i-1} - \tilde{b}_i'Q_i'{p}_{i-1})/(\bar{p}_{i-1}'\bar{p}_{i-1}). $$

   \item Compute the projection vector $p_i$ of $x$ over $ran([p_{i-1},A_i'])$:
       $$ p_i = x_i + \tilde{x}_i$$ and the inner product between $p_i$ and $x$:
       $$ c_i = \tilde{b}_i'\tilde{b}_i + \beta (c_{i-1} - \tilde{b}_iQ_i'p_{i-1})$$
   \item Go to next i.
  \end{list}

\item Output $p(=p_k)$ and $c(=c_k)$.
\end{list}

\end{alg}

 \subsection{Convergence analysis of SAP}

We now turn to the convergence of SAP (Algorithm \ref{alg:SAP}).  We first have the following conclusion.
\begin{thm} \label{thm:monotone}
Let $ \{x_{s}\}_1^t$ be any approximating vector sequence produced by the SAP process before convergence reached.
Then there exists \be \label{xi_norm}
   ||x_1 || < ||x_2|| < \cdots < || x_t ||  <  ||\bar{x}||
\ee
\end{thm}
\emph{Proof.}    \\
  Let $ \{p^s_j\}_{j=1}^k $ be the projection vector in the $s$-th AP iteration in SAP algorithm. Then we have
$x_i = p^i_k$ ($i=1,2,\cdots, \infty$).
Since $x_i \neq \bar{x}$ for any $i<t$, from Lemma \ref{lemma:AP4} we can always find some integer $j$ such that $A_j x_i \neq b_{j}$. Assuming $j$ is the first of this kind, then start from the $s+1$-th AP iteration in SAP algorithm, we have $  p^{i+1}_l= x_i $ $ (l=1,2,\cdots, j-1)$  and by Lemma \ref{lemma:AP1} we have  $$ ||x_i|| < ||p^{i+1}_j || $$
Note that we always have $ ||p^{i+1}_j || \le  ||p^{i+1}_k || = ||x_{i+1}||$, thus
 $$        ||x_i|| < ||p^{i+1}_j || \le || x_{i+1}||, \quad (i=1,2, \cdots, t) $$
and since $\bar{x}$ is the projection of $x$ over $ran(A')$ while $x_i$ is the projection of some subspace
$W_k $ of $ran(A')$, by Lemma \ref{lemma:3-2} we always have  $||x_i || < ||\bar{x}||$ for any $i$ $(1\le i \le t)$.
These complete the proof.  \hfill$\Box$

In order to prove the convergence of SAP method, we need the following conclusion.
\begin{lemma} \label{lemma:rr}
   Let $\{y_s\}_1^\infty$ be any convergent subsequence of $\{x_s\}_1^\infty$, and suppose
   $ \lim_{s\rightarrow \infty} y_s = y$.
   Then $y = \bar{x}$, where $\bar{x}$ is the projection of $x$ satisfying $Ax =b$.
\end{lemma}
\emph{Proof}. \\
  Since $y_s= x_{t_s}$ for some integer $t_s$ and $x_s = p^s_k$ ($s = 1,2,\cdots$),
  by Lemma \ref{lemma:AP1} we have
 $$ \barr{rl}
           ||y_{s+1}||^2 - ||y_s||^2 &= ||p^{t_{s+1}}_k ||^2 - || p^{t_s}_k||^2 \\
                                     &= ||p^{t_{s+1}}_k ||^2-||p^{t_s +1}_0||^2\\
                                     &>  || p^{t_s +1}_1 -  p^{t_s +1}_0||^2 \\
                                     &=  || p^{t_s +1}_1||^2 - ||p^{t_s +1}_0||^2 \\
                                     & = \tilde{x}' G \tilde{x} \earr
  $$
  where $ \tilde{x} = \bar{x} - x'u_0u_0$, $u_0 = p^{t_s}_0/||p^{t_s}_0||= y_s/||y_s||$ and $G = (\bar{A}_1\bar{A}_1')^{-1}$
  with $\bar{A}_1 = A_1(I - u_0u_0')$. i.e., we have
  \be\label{eq:mm}
       ||y_{s+1}||^2 - ||y_s||^2    >  \tilde{x}' G \tilde{x}
  \ee
  If $y \neq \bar{x}$, taking the limits on both sides of (\ref{eq:mm}) for $s $ approaching to infinity we have
     $$
       0 = \lim_{s \rightarrow \infty}    ||y_{s+1}||^2 - ||y_s||^2
       >  \lim_{s \rightarrow \infty} \tilde{x}'G \tilde{x}  >0
       $$
       since $G$ is symmetric positive definite and $\tilde{x}\neq 0$.  This is a contradiction and thus we must
       have $y = \bar{x}$.    \hfill$\Box$

\begin{thm}   \label{thm:convergence}
Let $ \{x_{s}\}_1^\infty$ be the approximating vector sequence produced by the SAP process.
Then
$$
     \lim_{s \rightarrow \infty} x_s = \bar{x}
$$
where
$\bar{x}$ stands for the projection of $x$ onto subspace $ran(A')$ of $R^n$. Particularly if $m =n$ and $A$ is
nonsingular, then we have $ \bar{x} = x$.
\end{thm}

\emph{Proof.   }
From Theorem \ref{thm:monotone} we see that sequence  $\{x_s\}_1^\infty$ is a bounded sequence.
If $ \displaystyle\lim_{s\rightarrow \infty} {x_s} $ does not exist, then it has at least two different cluster points $z_1$ and $z_2$
such that there are two subsequences of  $\{x_s\}_1^\infty$ approaches to $z_1$ and $z_2$. However by Lemma \ref{lemma:rr} we have $z_1 = \bar{x}$ and $z_2 = \bar{x}$. This is a contradiction. Hence we must have
      $$ \lim_{s\rightarrow \infty} {x_s} = \bar{x}. $$
Proof is completed. \hfill$ \Box$

\section{Some Acceleration Strategies}
We have observed from the preceding section that the convergence
speed of the simple iterative algorithm may not be very satisfactory
in general. In this section we are to design some accelerative
approaches for the SAP algorithm.

\subsection{Increase the Block Size}\label{sec:3_1}
An apparent approach is to simply increase the size of each block. The following table (Table \ref{table:iters2}) shows the iteration numbers needed for a convergent solution when $A = diag(-1,2,-1) \in R^{100\times 100}$ and the tolerance is set at $10^{-5}$ for the relative residual error. One can see that the number of iterations may drastically decrease when the size of blocks is slightly increased.  Unlike GMRES(m) with restarting where $m$ stands for the inner iteration numbers for each outer iteration and $m$ has to be very small comparing to the size of systems, this approach is viable since the size of each block can be selected much larger ( in case of sparse systems, one can choose the size of each block as large as $O(10^5)$ when LGO-based QS decomposition method \cite{pengDDM2009} is used to orthogonalize the block submatrices.

\begin{table}[h]\begin{center}
\caption{block-SAP-- iteration numbers needed for convergence}
\label{table:iters2}
\begin{tabular}{|c|c|c|c|c|c|c|c|c|c|c|}
\hline  block size  & $ 10$ &$ 15$ &$ 20 $ & $25 $ & $30$ & $35$ & $40$ &  $50$  \\ \hline
   iter\# &  $11404 $& $ 2994$ & $ 1020$& $ 443 $& $ 222$ &  $104 $& $ 57$   &$ 27$    \\ \hline
\end{tabular}
\end{center}
\end{table}

\subsection{A Modified SAP Approach}
Another option for accelerating the convergence is to add one simple step at the beginning of each loop in Algorithm 2 step 2. Specifically instead of using $x_{k}$ as the initial approximate solution to start another AP process, we first get the projection vector $p$ of $x$ onto the subspace $W = span\{ x_{k-1}, x_{k}\}$ as well as $c=x'p$ and then use it to replace $x_{k}$ and $c_{k}$ in Algorithm \ref{alg:SAP}. The details come as follows.

\begin{alg}\label{alg:MSAP} (A Modified Stationary Accumulated Projection Method-MSAP version 1). Let $A \in R^{n\times n}
$, $b \in R^n$. $\epsilon$ be a given tolerance. The following procedure produces an approximation
$p$ to the solution $x$ of system $Ax = b$.

\newcounter{msap1Counter1}
\begin{list}{Step\arabic{msap1Counter1}.} 
   { 
     \usecounter{msap1Counter1} 
     \setlength{\leftmargin}{3em}     
     \setlength{\labelwidth}{2. em}
     \setlength{\parsep}{0ex}         
     \setlength{\topsep}{1ex}         
     \setlength{\itemsep}{0.5ex}        
     \setlength{\labelsep}{.3em}     
     \setlength{\itemindent}{0em}    
     \setlength{\listparindent}{1em} 
   }
    \item Initialize $s$ as $s=0$, vector $x_0$ as $x_0 = \alpha A'b$, \mbox{$c_0 = \alpha b'b$}, $ t = ||b - Ax_0||/ ||b||$ where $\alpha = ||b||^2/||A'b||^2.$ Let $p = x_0$ and $c=c_0$.
    \item While $ t > \epsilon$
         \newcounter{msap1counter2}
         \begin{list}{Step\arabic{msap1counter2}.} 
         { 
     \usecounter{msap1counter2} 
     \setlength{\leftmargin}{2em}     
     \setlength{\labelwidth}{2. em}
     \setlength{\parsep}{0ex}         
     \setlength{\topsep}{1ex}         
     \setlength{\itemsep}{0.5ex}        
     \setlength{\labelsep}{.3em}     
     \setlength{\itemindent}{0em}    
     \setlength{\listparindent}{1em} 
      }

        \item Use Algorithm \ref{alg:APnew} to get a projection vector $x_{s+1}$ of $x$ and $c_{s+1}(=x'x_{s+1})$ with $p$ and $c$ the initial projection vector and the corresponding scalar taken as $x_{s}$ and $c_{s}$ respectively;
        \item[step 2.2] Calculate the projection $p$ of $x$ onto subspace $W = span\{ x_{s}, x_{s-1}\}$ and scalar  $c =x'p$. Rename $p$ as $x_{s+1}$ and $c_{s+1}$ respectively;
        \item[step 2.3] Calculate $t = ||b - Ax_{s+1}||/ ||b||;$
        \item[step 2.4]  $s = s+1$;
         \end{list}
    \item Output $p=(x_{s})$ and $c=(c_{s})$.
    \end{list}
\end{alg}
The following table (Table \ref{table:iters3}) shows the needed iteration numbers when running the same
example in subsection \ref{sec:3_1}.  It is easy to see that this simple acceleration technique works very well (comparing with the numbers in \mbox{Table \ref{table:iters2}}).

\begin{table}[h]\begin{center}
\caption{MSAP--iteration numbers needed for convergence}
\label{table:iters3}
\begin{tabular}{|c|c|c|c|c|c|c|c|c|c|c|}
\hline  block size  & $ 10$ &$ 15$ &$ 20 $ & $25 $ & $30$ & $35$ & $40$ &  $50$  \\ \hline
   iter\# &    $2134$ & $403$ & $134$ & $69$ & $38$&  $34$ & $18$ & $15 $  \\ \hline
\end{tabular}
\end{center}
\end{table}

It is attempting to increase the dimension of the subspace $W$ in Algorithm \ref{alg:MSAP} step 2.2 to get a better
convergence speed. Unfortunately this seems not work since the ``distance" between $\{x_{s}\}_1^m$ are not far enough and thus the submatrix formed by these vectors tends to be very ill-conditioned, which eventually makes the idea not work well. The following subsection depicts an alternative option for this idea.

\subsection{ A Varying Subspace Method}

It is attractive to use a subspace $W$ with larger dimension than that of subspace $W$ in Algorithm  \ref{alg:MSAP} step 2.2 to develop an accelerative method for SAP. However we have noted that as iteration goes on, the matrix formed by the successive SAP projection tends to be ill-conditioned. Hence we plan to use a more flexible strategy to handle the ill-conditioned systems. Our intension is to use a detector to check the conditioning of an intermediate matrix $H$, and then arrange the dimension of the subspace $W$ accordingly. The details are described in the following algorithm.

\begin{alg}\label{alg:MSAP2} (A Modified Stationary Accumulated Projection Method-MSAP version 2). Let $A \in
R^{n\times n}
$, $b \in R^n$. $\epsilon$ be a given tolerance. The following procedure produces an approximation
$p$ to the solution $x$ of system $Ax = b$.
\newcounter{numcount1}
\begin{list}{Step\arabic{numcount1}.} 
   { 
     \usecounter{numcount1} 
     \setlength{\leftmargin}{3em}     
     \setlength{\labelwidth}{2. em}
     \setlength{\parsep}{0ex}         
     \setlength{\topsep}{1ex}         
     \setlength{\itemsep}{0.5ex}        
     \setlength{\labelsep}{.3em}     
     \setlength{\itemindent}{0em}    
     \setlength{\listparindent}{1em} 
   }
\item Initialize $s$ as $s=0$, vector $x_{0}$ as $x_{0} = \alpha A'b$, \mbox{$c_{0} = \alpha b'b$}, $ t = ||b - Ax_{0}||/ ||b||$ where $ \alpha = ||b||^2/||A'b||^2$. Let $p = x_{0}$ and $c=c_{0}$,
     and $m$ be a small predetermined integer.
\item While $ t > \epsilon$
        \newcounter{newcounter2} \begin{list}{Step 2.\arabic{newcounter2}.} 
           { 
             \usecounter{newcounter2} 
             \setlength{\leftmargin}{2.em}     
             \setlength{\labelwidth}{2em}
             \setlength{\parsep}{0ex}         
             \setlength{\topsep}{1ex}         
             \setlength{\itemsep}{0.5ex}        
             \setlength{\labelsep}{0.3em}     
             \setlength{\itemindent}{.5em}    
             \setlength{\listparindent}{1.4em} 
           }
        \item Use Algorithm \ref{alg:APnew} to get a projection vector $pn$ of $x$ and $cn (=x'x_{s+1})$ with $p$ and $c$ as the initial projection vector and the corresponding scalar.
        \item Store $pn$ as a row vector in matrix $H$ and $cn$ into a column vector $L$.
        \item If $H$ contains $m$ row vectors\hfill

            \newcounter{newcounter3}
                \begin{list}{Step 2.\arabic{newcounter2}.\arabic{newcounter3}} 
                   { 
                     \usecounter{newcounter3} 
                     \setlength{\leftmargin}{2.0em}     
                     \setlength{\labelwidth}{2. em}
                     \setlength{\parsep}{0ex}         
                     \setlength{\topsep}{1ex}         
                     \setlength{\itemsep}{0.5ex}        
                     \setlength{\labelsep}{0.3em}     
                     \setlength{\itemindent}{.5em}    
                     \setlength{\listparindent}{1.4em} 
                   }
              \item if $H$ is well-conditioned,\hfill
                    \newcounter{newcounter4}
                    \begin{list}{Step 2.3.1.\arabic{newcounter4}}
                    { 
                      \usecounter{newcounter4} 
                      \setlength{\leftmargin}{2.em}     
                       \setlength{\labelwidth}{2. em}
                      \setlength{\parsep}{0ex}         
                      \setlength{\topsep}{1ex}         
                      \setlength{\itemsep}{0.5ex}        
                      \setlength{\labelsep}{0.3em}     
                      \setlength{\itemindent}{.5em}    
                      \setlength{\listparindent}{2.5em} 
                    }
                      \item update $x_{s+1}$ as the projection of $x$ on subspace $W = ran(H)$ and the scalar $c_{s+1} = x' x_{s+1}$
                      \item remove the first row vector of $H$ and the first element in vector $L$ correspondingly;
                  \end{list}
              \item[step 2.3.2] if $H$ is ill-conditioned\hfill
                  \newcounter{newcounter5}
                    \begin{list}{Step 2.3.2.\arabic{newcounter5}}
                    { 
                      \usecounter{newcounter5} 
                      \setlength{\leftmargin}{2.5em}     
                      \setlength{\parsep}{0ex}         
                      \setlength{\topsep}{1ex}         
                      \setlength{\itemsep}{0.5ex}        
                      \setlength{\itemindent}{.5em}    
                    }
                    \item update $x_{k+1}$ as the projection of $x$ on subspace
                        $W = span(x_{s},pn)$  and the scalar $c_{s+1} = x' x_{s+1}$.
                    \item remove all but the first row vectors of $H$ and the
                        elements in vector $L$ correspondingly.
                    \end{list}
            \end{list}
        \item[step 2.4] if $H$ contains less than $m$ row vectors,
             update $x_{s+1}$ as the projection of $x$ on subspace $W = span(x_{s},pn)$  and the scalar $c_{s+1} = x' x_{s+1}.$
        \item[step 2.5] Calculate $t = ||b - Ax_s||/ ||b||;$
        \item[step 2.6]  Set $s = s+1$, $p = x_{s+1}$ and $c = c_{s+1}$.
        \end{list}
\item Output $p(=x_{s})$ and $c_{s}$.
\end{list}

\end{alg}

The following table (Table \ref{table:iters4}) shows the astonishing acceleration speed of convergence when we use Algorithm \ref{alg:MSAP2} to solve the same problem in the preceding subsection.

\begin{table}[h]\begin{center}
\caption{MSAP2--iteration numbers needed for convergence}
\label{table:iters4}
\begin{tabular}{|c|c|c|c|c|c|c|c|c|c|c|}
\hline  block size  & $ 10$ &$ 15$ &$ 20 $ & $25 $ & $30$ & $35$ & $40$ &  $50$  \\ \hline
   iter\# &    $185$ & $102$ & $42$ & $30$ & $16$&  $14$ & $10$ & $7 $  \\ \hline
\end{tabular}
\end{center}
\end{table}

\section{Numerical Experiments}
In this section we will show some application of the aforementioned
SAP methods and we compare the results with those produced by GMRES---a benchmark Krylov subspace method. We use MSAP and GMRES to calculate the systems.

As the first example, we use the two-point boundary value problem
\be \left\{
\barr{l}
    (a(t)u^\prime(t))^\prime + b(t)u(t) = f(t), \quad t \in (0,1)  \\
    u(0) = u(1) = 0.
\earr
\right.
\ee
This equation represents some important practical problems such as chord balancing, elastic beam problems, etc.
We use finite element method to get the numerical solution to the system, which ends up with a linear system of equations in the form of (\ref{eq:1}) with $n$ unknowns, where $n$ stands for the number of grids which divide the interval $(0,1)$ into
$n+1$ equal-sized subintervals $(x_i,x_{i+1})$ ($i = 1,2,\cdots n$). We use the linear interpolation function at each grid point to construct the finite element space $V_0^h$, test functions are also from  $V_0^h$.

In our test we take $a(t) = 1+t$, $b(t) = t$ and $f(t)$ is taken so that the exact solution to the system is
$ u(t) = t(1-t)e^{2+t}$.  By using the aforementioned finite element space $V_0^h$, we get a linear system of equation $A x = b$ with $A$ as a symmetric tridiagonal matrix in $R^{n \times n}$, where $n$ is taken as $200$.
The block size for the MSAP method (version 2) is set to be the square root of the restart number $m$ of GMRES(m) multiplied by the number of unknowns $n$ so that submatrices of $A$ have roughly the same number of non-zero elements as those in Krylov subspace matrices formed in GMRES process.

Table \ref{table:iters5}  shows the comparison of iteration numbers needed for convergence, the relative error in terms of $||x - approx.x||/||x||$ (where $approx.x$ stands for the approximate solutions obtained by using GMRES and MSAP2 respectively) and  the relative residual error in terms of \mbox{$||b - A* approx.x||/||b||$ }
 is used for the convergence criteria. We also observed that the time cost in these example also show some advantage of MSAP over GMRES as shown in Table \ref{table:iters5}.

\begin{table}[ht]\begin{center}
\caption{Comparison between MSAP and GMRES}
\label{table:iters5}
\begin{tabular}{|c|c|c|c|c|c|c|c|c|}
\hline
\multicolumn{2}{|c|}{settings} & \multicolumn{2}{|c|}{iter. \#}  &  \multicolumn{2}{|c|}{time(in s)} &  \multicolumn{2}{|c|}{rel. error} \\ \hline
msap &gmres & msap & gmres & msap & gmres & msap & gmres \\\hline
 blk\_size& restart  & m=5   &(out,in)  &  &   &   &  \\ \hline
20&2 & 200 & (2000,2)  &   1.2690 & 0.6700  &  7.02e-7  &  6.96e-5   \\ \hline
30&5 & 200 &(2000,5)  &   0.6340&   1.6430 & 3.49e-7  &  5.68e-5\\ \hline
40&8 & 50 &(1415,7 )   &   0.1320&  1.8180  & 9.57e-7 &   4.92e-5\\ \hline
50&13 & 33 & (538,3 )  &   0.0800& 1.3780  & 3.32e-7  &   4.4e-5\\ \hline
60&18 &22 &(282,18)  &   0.0510& 1.1880  &  2.41e-7  &  4.01e-5\\ \hline
70&25 &17 &(148,24 )  &  0.0380& 1.1750   & 5.06e-7  &  3.71e-5\\ \hline
80&32  &13&(91,28 ) &  0.0380&  0.9950  &3.01e-8 &   3.48e-5\\ \hline
\end{tabular}
\end{center}
\end{table}

From the above table it seems that MSAP has a better relative error level than that of GMRES
at the same relative residual level. We have to point out here that the construction of this
test system is made so that the solution has rich eigenvector components corresponding to the
smallest eigenvalues of the coefficient matrix $A$. In case the condition number $cond(A)$ is relatively
small (say, less than $O(10^3)$), GMRES  outperforms MSAP in terms of time costs and flops in our tests, while in case
the condition number of the coefficient matrices are larger than $O(10^3)$, MSAP generally outperforms GMRES
in most of our test cases.

Table $7$\, is a comparison between MSAP and block Jacobi method for a system with coefficient matrix $ A \in R^{200 \times 200}$. The block size for both methods are chosen as exactly the same. One can see a much less iteration number needed for MSAP than that of block Jacobi method in each case of block size, while the relative
errors obtained by MSAP are much better than those obtained by block Jacobi method.

\begin{table}[h]
\begin{center}
\caption{ iteration numbers needed for convergence, tol =$10^{-5}$}
\begin{tabular}{|c|c|c|c|c|c|c|c|c|c|} \hline
 blk\_size & \multicolumn{2}{|c|}{time(in s)} &  \multicolumn{2}{|c|}{iteration\#} &   \multicolumn{2}{|c|}{rel. error}  \\\hline
                   &   Jacobi &  msap        &  Jacobi &  msap        & Jacobi &  msap        \\\hline
            10    &    5.666  &     12.166   &      7836   &      1745&  6.9553e-005 & 7.0191e-007    \\\hline
            15    &    2.769  &       3.68   &      5347   &       830 & 5.6761e-005 & 3.4931e-007     \\\hline
            20    &    1.547   &     1.489   &      4082   &       390 &  4.921e-005 & 9.5735e-007      \\\hline
            25    &    1.018  &      0.587   &      3316   &       185 & 4.3997e-005&  3.3189e-007    \\\hline
            30    &    0.781  &      0.379   &      2806   &       130 & 4.0135e-005&  2.4148e-007    \\\hline
            35    &    0.589  &      0.222  &       2440   &        85 & 3.7142e-005&  5.0644e-007     \\\hline
            40    &    0.446  &      0.134  &       2159   &        55 & 3.4765e-005 & 3.0083e-008     \\\hline
            45   &     0.408  &      0.102  &       1946    &       45 & 3.2554e-005 &  3.132e-008     \\\hline
 \end{tabular}
\end{center}
 \end{table}

\section{Comments and Summary}

In this paper we present a new type of iterative methods for solving linear system of \mbox{equations}. This might be the first type of methods which do not belong to the category of extended Krylov subspace methods as we mentioned above. It can overcome some shortcomings of Krylov subspace methods and exhibit better performance in our test problems. We need to mention that convergence speed of these algorithms deteriorates when the number of subdividing blocks of the coefficient matrix exceeds 20, a remedy to this is to embed
an inner loop in the AP process, which will cause more flops but the obtained time efficiency payoff
this costs in our tests.
 A  relevant issue is the study of detailed error analysis for the SAP algorithm which may leads to a deep insight error estimation for each AP process in terms of subspace distance.
We need to mention here that the SAP algorithm is nothing but a ``horizontal" application of the AP process, i.e., the AP processes are always applied to the original linear system instead of residual equations. We find that a vertical application of AP process is  also possible and the results are to appear in our later work.

\bibliographystyle{plain}

\bibliography{peng_ref}
\end{document}